\documentclass{article}
\usepackage[a4paper]{geometry}
\usepackage{amsmath,amssymb,amsthm,enumitem,environ,color,mathrsfs,tikz}
\usepackage[initials]{amsrefs}
    
\newtheorem{theorem}{Theorem}[section]
\newtheorem{lemma}[theorem]{Lemma}

\theoremstyle{definition}

\theoremstyle{plain}

\numberwithin{equation}{section}
\numberwithin{figure}{section}

\renewcommand{\geq}{\geqslant}
\renewcommand{\leq}{\leqslant}

\setlist[enumerate]{leftmargin=20pt,itemsep=0pt,topsep=0pt}
\setlist[enumerate,1]{label=\emph{(\roman*)},ref={(\roman*)}}
\setlist[enumerate,2]{label=\emph{(\alph*)},ref={(\alph*)}}

\makeatletter
\renewcommand\section{\@startsection {section}{1}{\z@}%
                                   {-3.5ex \@plus -1ex \@minus -.2ex}%
                                   {1.3ex \@plus.2ex}%
                                   {\normalfont\large\scshape}}
\makeatother

\usetikzlibrary{calc}
\usetikzlibrary{matrix}
\usetikzlibrary{decorations.markings}
\usetikzlibrary{arrows.meta}
\usetikzlibrary{arrows,shapes,positioning}
\tikzstyle arrowstyle=[scale=1.3]
\tikzstyle directed=[postaction={decorate,decoration={markings,
  mark=at position 0.6 with {\arrow[arrowstyle]{{latex}}}}}]
\tikzstyle directedplus=[postaction={decorate,decoration={markings,
  mark=at position 0.7 with {\arrow[arrowstyle]{{latex}}}}}]

\tikzset{pointer/.style 2 args={draw,fill,single arrow,
    single arrow tip angle=45,
    single arrow head extend=#1,
    single arrow head indent=0pt,
    inner sep=0pt,
    rotate=#2}}

\definecolor{grey}{RGB}{170,170,170}
\definecolor{gammacol}{RGB}{33,120,33}
\definecolor{deltacol}{RGB}{44,90,160}
\definecolor{alphacol}{RGB}{170,0,0}
\definecolor{betacol}{RGB}{171,55,200}

\newcount\recurdepth
\def\fareyrecur[#1]#2#3#4#5#6#7{
  \recurdepth=#6
  \ifnum\the\recurdepth>1\relax
    \advance\recurdepth by-1\relax
    \edef\tempnum{\number\numexpr#2+#4\relax}
    \edef\tempden{\number\numexpr#3+#5\relax}
    \pgfmathparse{\tempnum/\tempden}\edef\temp{\pgfmathresult}

    \ifnum #7=0 
      \ifnum\tempnum>0
        \node[below=1pt,scale=1]at({(\temp)},0){$\frac{\tempnum}{\tempden}$};
      \else 
        \ifnum\the\recurdepth>1\relax 
          \edef\abstempnum{\number\numexpr -\tempnum\relax}
          \node[below=1pt,scale=1]at({(\temp)},0){$-\frac{\abstempnum}{\tempden}\phantom{-}$};
        \fi 
      \fi 
      \draw[#1] ({(\temp)},0) arc (180:0:{((#4/#5)-\temp)*0.5});
      \draw[#1] ({(\temp)},0) arc (0:180:{(\temp-(#2/#3))*0.5});
    \fi

    \begingroup
      \edef\ttempup{\noexpand\fareyrecur[#1]{\tempnum}{\tempden}{#4}{#5}{\the\recurdepth}{#7}}
      \edef\ttempdown{\noexpand\fareyrecur[#1]{#2}{#3}{\tempnum}{\tempden}{\the\recurdepth}{#7}}
      \ttempup\ttempdown
    \endgroup
  \fi
}

\def\fareygraph[#1]#2#3#4{
  \draw[#1] ({#2},0) -- ({#3+1},0);
  \foreach \n in {#2,...,#3}{
    \draw[#1] ({\n},0.75) -- ({\n},0);
    \ifnum #4=0 
      \node[below=1pt,scale=1,black] at ({\n},0) {$\frac{\n}{1}$};	
    \fi
    \draw[#1] ({(\n)},0) arc (180:0:0.5);
    \fareyrecur[#1]{\n}{1}{\n+1}{1}{4}{#4}	
  }
  \edef\nplusone{\number\numexpr#3+1\relax}
  \draw[#1] ({\nplusone},0.75) -- ({\nplusone},0) node[below=4pt,scale=1] {$\frac{\nplusone}{1}$};	
}

\title{ \vspace{-5ex}\bf \Large Necessary and sufficient conditions for convergence of integer continued fractions}

\makeatletter
\renewcommand*{\@fnsymbol}[1]{\hspace*{-10pt}}
\makeatother

\author{Ian Short and Margaret Stanier\thanks{School of Mathematics and Statistics, The Open University, Milton Keynes, MK7 6AA, United Kingdom.}}

\date{\vspace{-5ex}}

\begin{document}

\maketitle

\begin{abstract}
Fundamental to the theory of continued fractions is the fact that every infinite continued fraction with positive integer coefficients converges; however, it is unknown precisely which continued fractions with integer coefficients (not necessarily positive) converge. Here we present a simple test that determines whether an integer continued fraction converges or diverges. In addition, for convergent continued fractions the test specifies whether the limit is rational or irrational. 

An attractive way to visualise integer continued fractions is to model them as paths on the Farey graph, which is a graph embedded in the hyperbolic plane that induces a tessellation of the hyperbolic plane by ideal triangles. With this geometric representation of continued fractions our test for convergence can be interpreted in a particularly elegant manner, giving deeper insight into the nature of continued fraction convergence.
\end{abstract}

\section{Introduction}

It is well known that every infinite continued fraction 
\[
(b_0,b_1,\dotsc) = b_0+ \cfrac{1}{b_1+\cfrac{1}{b_2+\cfrac{1}{\raisebox{-1ex}{$b_3+\dotsb$}}}}
\]
with positive integer coefficients converges (to an irrational number). However, if we stipulate that the coefficients are integers, but not necessarily positive integers, then the continued fraction need not converge; for example, the periodic continued fraction $(1,-1,1,-1,\dotsc)$ diverges -- its sequence of convergents oscillates between the three values $1$, $0$ and $\infty$. This paper gives a simple set of necessary and sufficient conditions for an integer continued fraction to converge. Furthermore, for each continued fraction that does converge, the conditions specify whether the value of the continued fraction is rational or irrational.

There is an extensive literature on the convergence of continued fractions -- from Perron's classic  treatise \cite{Pe1954,Pe1957} to more recent works such as \cite{BoMc2007,Lo2016,LoWa2008} -- and there are a wealth of algorithms for generating different types of integer continued fractions (see, for example, \cite[Chapter 4]{IoKr2002} for algorithms in the metric theory of continued fractions, \cite{Ka2013,KaUg2007} for algorithms in the geometric theory of continued fractions, and \cite{Da2015,DaNo2014} for a general class of algorithms that yield continued fractions with integer and Gaussian integer coefficients). A wide variety of tests for convergence of continued fractions with real and complex coefficients are known. For instance, if $|b_n|\geq 2$ for $n=0,1,\dotsc$, then $(b_0,b_1,\dotsc)$ converges -- even if the coefficients $b_n$ are complex numbers -- and some slightly stronger tests are known when the coefficients are integers (such as \cite[Lemma~1.1]{KaUg2005}). Here we present the first test that determines precisely whether or not any integer continued fraction converges.

It is more convenient to state all our results in terms of negative infinite continued fractions
\[
[b_0,b_1,\dotsc] = b_0- \cfrac{1}{b_1-\cfrac{1}{b_2-\cfrac{1}{\raisebox{-1ex}{$b_3-\dotsb$}}}}
\]
with integer coefficients. We can easily switch between regular continued fractions and negative continued fractions using the formula
\[
(b_0,b_1,b_2,b_3,\dotsc) = [b_0,-b_1,b_2,-b_3,\dotsc].
\]
These two continued fractions have the same sequence of convergents, so one converges if and only if the other does, and if they do converge then they converge to the same value.

Henceforth we use the phrase `continued fraction' to refer to a negative continued fraction with integer coefficients. Usually our continued fractions will be infinite, but on occasion we will need finite continued fractions $[b_0,b_1,\dots,b_n]$, with the obvious use of notation. 

A continued fraction $[b_0,b_1,\dotsc]$ is considered to converge if its sequence of convergents $v_n=[b_0,b_1,\dots,b_n]$, for $n=0,1,\dotsc$, converges in the extended real line $\mathbb{R}_\infty$ (where $\mathbb{R}_\infty=\mathbb{R}\cup\{\infty\}$). It diverges if it does not converge. The limit, when it exists, is called the value of the continued fraction. The convergents belong to the set of extended rationals $\mathbb{Q}_\infty$ (where $\mathbb{Q}_\infty=\mathbb{Q}\cup\{\infty\}$), but the value of the continued fraction lies in $\mathbb{R}_\infty$, and it could be irrational.

We say that two continued fractions  $[b_0,b_1,\dotsc]$ and  $[c_0,c_1,\dotsc]$ have the same \emph{tail} if there are positive integers $r$ and $s$ such that the two sequences $b_r,b_{r+1},\dotsc$ and $c_s,c_{s+1},\dotsc$ coincide.

Let $\mathscr{C}$ denote the collection of all continued fractions $[b_0,b_1,\dotsc]$. We define a function $\Phi\colon\mathscr{C}\longrightarrow\mathscr{C}$ as follows. If $b_n\neq 0,1,-1$ for every \emph{positive} integer $n$ (ignore $b_0$), then $\Phi$ fixes $[b_0,b_1,\dotsc]$. Otherwise, let $m$ be the least positive integer for which $b_m$ is $0$, $1$ or $-1$. Then
\[
\Phi([b_0,b_1,\dotsc])=
\begin{cases}
[b_0,b_1,\dots,b_{m-2},b_{m-1}+b_{m+1},b_{m+2},b_{m+3},\dotsc], & \text{if $b_m=0$,}\\
[b_0,b_1,\dots,b_{m-2},b_{m-1}-1,b_{m+1}-1,b_{m+2},b_{m+3},\dotsc], & \text{if $b_m=1$,}\\
[b_0,b_1,\dots,b_{m-2},b_{m-1}+1,b_{m+1}+1,b_{m+2},b_{m+3},\dotsc], & \text{if $b_m=-1$.}
\end{cases}
\]
In each case, $\Phi$ removes the coefficient $b_m$, and adjusts the two coefficients on either side, merging them when $b_m=0$. The operations induced by $\Phi$ are familiar in the theory of continued fractions, where they are sometimes referred to as `singularization' operations (see, for example \cite[Section~4.2]{IoKr2002}). Towards the end of this introduction we give a geometric description of the function $\Phi$ which will help to explain its definition.

Given $[b_0,b_1,\dotsc]$ in $\mathscr{C}$, we define $p^{(n)}$ to be the first positive integer position (not zero) at which a $0$, $1$ or $-1$ appears in the continued fraction $\Phi^n([b_0,b_1,\dotsc])$, and we define $p^{(n)}$ to be $\infty$ if there is no such position. Let $p=\liminf p^{(n)}$. Suppose for the moment that $p=\infty$. In this case, for each nonnegative integer $k$, the sequence of integers obtained by taking the $k$th coefficient of $\Phi^n([b_0,b_1,\dotsc])$, for $n=1,2,\dotsc$, eventually fixes on some value $b_k^*$. Thus we obtain a limit continued fraction $[b_0^*,b_1^*,\dotsc]$, where $|b^*_n|\geq 2$, for $n=1,2,\dotsc$. Suppose now that $p<\infty$. In this case, we define $q^{(n)}$ to be the modulus of the coefficient of $\Phi^n([b_0,b_1,\dotsc])$ at position $p-1$, for $n=1,2,\dotsc$, to give a sequence of nonnegative integers $(q^{(n)})$.

We can now state our principal theorem, which uses the notation developed in the preceding paragraph.

\begin{theorem}\label{thm1}
Let $[b_0,b_1,\dotsc]$ be a continued fraction.
\begin{enumerate}
\item Suppose that $p^{(n)}\rightarrow \infty$.
\begin{enumerate}
\item If $[b^*_0,b^*_1,\dotsc]$ has the same tail as $[2,2,\dotsc]$ or $[-2,-2,\dotsc]$, then $[b_0,b_1,\dotsc]$ converges to a rational.
\item Otherwise, $[b_0,b_1,\dotsc]$ converges to an irrational.
\end{enumerate}
\item Suppose that $p^{(n)}\not\rightarrow \infty$.
\begin{enumerate}
\item If $q^{(n)}\rightarrow\infty$, then $[b_0,b_1,\dotsc]$ converges to an extended rational.
\item If $q^{(n)}\not\rightarrow\infty$, then $[b_0,b_1,\dotsc]$ diverges.
\end{enumerate}

\end{enumerate}
\end{theorem}

It is worth highlighting a corollary of part (i) of Theorem~\ref{thm1} for continued fractions with no coefficients (other than perhaps $b_0$) equal to $0$, $1$ or $-1$.

\begin{theorem}\label{thm2}
Suppose that $|b_n|\geq 2$, for $n=1,2,\dotsc$. Then the  continued fraction $[b_0,b_1,\dotsc]$ converges to a rational if it has the same tail as $[2,2,\dotsc]$ or $[-2,-2,\dotsc]$, and otherwise it converges to an irrational.
\end{theorem}

This theorem generalizes the known result of the same type that has $b_n\geq 2$, for $n=1,2,\dotsc$ (see, for example, \cite[Section~1]{Ka1996}).

We provide some examples to illustrate how Theorem~\ref{thm1} can be applied. Consider first the continued fraction
\[
[b_0,b_1,\dotsc] = [3,1,3,\textcolor{grey}{4},\textcolor{grey}{1},\textcolor{grey}{2},\textcolor{grey}{3},5,1,2,2,3,\textcolor{grey}{6},\textcolor{grey}{1},\textcolor{grey}{2},\textcolor{grey}{2},\textcolor{grey}{2},\textcolor{grey}{3},\dotsc].
\]
Some numbers are shaded to indicate how the pattern of coefficients continues. Now, if these were the coefficients of a regular continued fraction, then that continued fraction would converge, because all the coefficients are positive. Remember, however, that $[b_0,b_1,\dotsc]$ is a negative continued fraction, and to determine whether it converges we can use Theorem~\ref{thm1}. To do this, we apply $\Phi$ repeatedly to obtain the sequence $\Phi^n([b_0,b_1,\dotsc])$, for $n=1,2,\dotsc$, which is
\[
[2,2,\textcolor{grey}{4},\textcolor{grey}{1},\textcolor{grey}{2},\textcolor{grey}{3},5,1,\dotsc] \rightarrow
[2,2,\textcolor{grey}{3},\textcolor{grey}{1},\textcolor{grey}{3},5,1,\dotsc] \rightarrow
[2,2,\textcolor{grey}{2},\textcolor{grey}{2},5,1,\dotsc] \rightarrow \dotsb.
\] 
Continuing in this way we see that $p^{(n)}\to\infty$ and the limit continued fraction is
\[
[b_0^*,b_1^*,\dotsc] = [2,2,\dotsc].
\]
Hence $[b_0,b_1,\dotsc]$ converges to a rational number (namely 1, the value of $[2,2,\dotsc]$).

Consider now the continued fraction
\[
[b_0,b_1,\dotsc] = [1,2,1,3,1,4,1,5,1,\dotsc].
\]
Applying $\Phi$ repeatedly yields the sequence 
\begin{align*}
&[1,1,2,1,4,1,5,1,6,1,7,\dotsc]\rightarrow
[0,1,1,4,1,5,1,6,1,7,\dotsc]\rightarrow
[-1,0,4,1,5,1,6,1,7,\dotsc] \\
&\rightarrow[3,1,5,1,6,1,7,\dotsc]\rightarrow
[2,4,1,6,1,7,\dotsc]\rightarrow
[2,3,5,1,7,\dotsc]\rightarrow \dotsb.
\end{align*}
Again we see that $p^{(n)}\to\infty$, and this time the limit continued fraction is
\[
[b_0^*,b_1^*,\dotsc] = [2,3,4,5,\dotsc].
\]
Hence $[b_0,b_1,\dotsc]$ converges to an irrational number. 

Next consider the continued fraction
\[
[b_0,b_1,\dotsc]= [1, 0, 2, 0, 3, 0, 4, 0,\dotsc].
\]
Applying $\Phi$ repeatedly gives
\[
[3, 0, 3, 0, 4, 0,5,0\dotsc] \rightarrow
[6, 0, 4, 0,5,0\dotsc] \rightarrow
[10, 0,5,0\dotsc] \rightarrow \dotsb.
\]
By continuing with this process it becomes clear that, for each positive integer $n$, the coefficient of $\Phi^n([b_0,b_1,\dotsc])$ with position~$1$ is $0$, and the coefficient with position~$0$ is  $\tfrac12(n+1)(n+2)$. Hence $p=\liminf p^{(n)}=1$ and $q^{(n)}\to\infty$. Therefore $[b_0,b_1,\dotsc]$ converges to an extended rational -- to infinity, in fact (this could easily be ascertained by other means).

Last consider the continued fraction
\[
[b_0,b_1,\dotsc]= [3,0,-3,\textcolor{grey}{3},\textcolor{grey}{3},\textcolor{grey}{0},\textcolor{grey}{-3},\textcolor{grey}{-3},3,3,3,0,-3,-3,-3,\dotsc].
\]
Once more, by applying $\Phi$ repeatedly we obtain
\begin{align*}
&[0,\textcolor{grey}{3},\textcolor{grey}{3},\textcolor{grey}{0},\textcolor{grey}{-3},\textcolor{grey}{-3},3,3,3,0,-3,-3,-3,\dotsc] \rightarrow
[0,\textcolor{grey}{3},\textcolor{grey}{0},\textcolor{grey}{-3},3,3,3,0,-3,-3,-3,\dotsc] \\
& \rightarrow
[0,\textcolor{grey}{0},3,3,3,0,-3,-3,-3,\dotsc] \rightarrow
[3,3,3,0,-3,-3,-3,\dotsc] \rightarrow \dotsb.
\end{align*}
In this case the sequence of coefficients of $\Phi^n([b_0,b_1,\dotsc])$ with position~$0$, for $n=1,2,\dotsc$, takes the value $0$ infinitely often, as does the sequence of coefficients with position~1. Hence $p=\liminf p^{(n)}=1$ and $q^{(n)}\not\rightarrow\infty$, so  $[b_0,b_1,\dotsc]$ diverges.

To prove Theorem~\ref{thm1} we make use of the Farey graph, which we now define. Let $\mathbb{H}$ be the (open) upper half-plane; this is a model of the hyperbolic plane when it is endowed with the Riemannian metric $|dz|/\text{Im}\, z$.  Each hyperbolic line in this model of the hyperbolic plane is either the upper half of a circle centred on the real axis or a half-line in $\mathbb{H}$ orthogonal to the real axis. We are interested in the set of extended rationals $\mathbb{Q}_\infty$, which lie on the boundary $\mathbb{R}_\infty$ of $\mathbb{H}$. We use the term \emph{reduced rational} to describe an expression $a/b$ in which $a$ and $b$ are coprime integers and $b>0$. For convenience, we consider the expression $1/0$ to be a reduced rational also, representing the point $\infty$.

\begin{figure}[ht]
\centering
\begin{tikzpicture}[scale=4]
\begin{scope}
    \clip(-0.9,-0.2) rectangle (2.3,0.75);
	\fareygraph[thin]{-1}{2}{0}
\end{scope}
\end{tikzpicture}
\caption{Part of the Farey graph}
\label{fig1}
\end{figure}
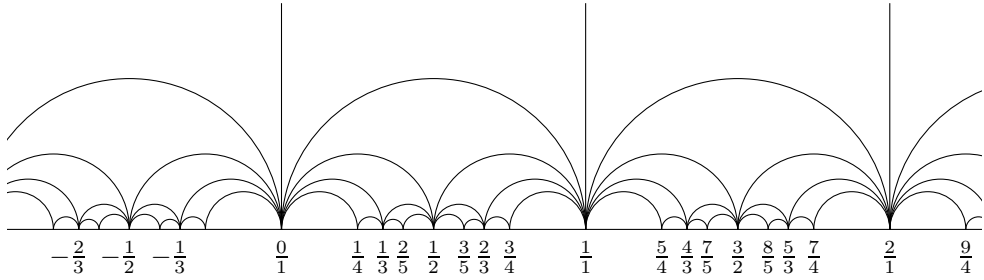

The \emph{Farey graph} $\mathscr{F}$ is the graph with vertices $\mathbb{Q}_\infty$ and edges comprising those pairs $a/b$ and $c/d$ of reduced rationals for which $ad-bc=\pm 1$. We represent the edge incident to $a/b$ and $c/d$ by the unique hyperbolic line in $\mathbb{H}$ between those two boundary points. The collection of all edges creates a tessellation of $\mathbb{H}$ by triangles, part of which is shown in Figure~\ref{fig1}.

The \emph{modular group} $\Gamma$ comprises all M\"obius transformations of the form
\[
f(z) = \frac{az+b}{cz+d}, 
\]
where $a,b,c,d\in\mathbb{Z}$ and $ad-bc=1$. The modular group acts on the extended complex plane in the usual way, with the standard conventions regarding the point $\infty$. It also acts on $\mathbb{R}_\infty$, on $\mathbb{Q}_\infty$, and on $\mathbb{H}$ as a group of orientation-preserving hyperbolic isometries. One can check that $f$ preserves adjacency between vertices of $\mathscr{F}$, so $\Gamma$ acts on $\mathscr{F}$ too; in fact, each element of $\Gamma$  induces a graph automorphism of $\mathscr{F}$.

Given two vertices $u$ and $v$ of $\mathscr{F}$, we write $u\sim v$ if $u$ and $v$ are adjacent. An \emph{infinite path} in $\mathscr{F}$ is a sequence $v_0,v_1,\dotsc$ of vertices of $\mathcal{F}$ such that $v_i\sim v_{i+1}$, for $i=0,1,\dotsc$. We denote this path by $\langle v_0,v_1,\dotsc\rangle$. Occasionally we consider finite paths in $\mathscr{F}$, which are defined in the obvious way, with similar notation. An infinite path $\langle v_0,v_1,\dotsc\rangle$ is said to \emph{converge} in $\mathscr{F}$ if the sequence $v_0,v_1,\dotsc$ converges in $\mathbb{R}_\infty$. Note that we allow paths to pass through the same vertex more than once; in graph theory the term `walk' is often used in place of `path'.

For a continued fraction $[b_0,b_1,\dotsc]$, we define $s_n(z)=b_n-1/z$ and $S_n=s_0\circ s_1\circ \dotsb \circ s_n$, for $n=0,1,\dotsc$. The maps $s_n$ and $S_n$ belong to the modular group $\Gamma$. Let us define $v_n=S_n(\infty)$, for $n=0,1,\dotsc$. Then $v_n$ is the $n$th convergent of $[b_0,b_1,\dotsc]$. 

Since $0\sim\infty$ and $S_n\in\Gamma$ it follows that $S_n(0)\sim S_n(\infty)$. But $S_n(0)=S_{n-1}\circ s_n(0)=S_{n-1}(\infty)$, so we can see that $\langle \infty, v_0,v_1,\dotsc\rangle$ is an infinite path in $\mathscr{F}$. Conversely, given an infinite path $\langle \infty, v_0,v_1,\dotsc\rangle$ in $\mathscr{F}$ there is a unique infinite continued fraction $[b_0,b_1,\dotsc]$ with sequence of convergents $v_0,v_1,\dotsc$. This is relatively straightforward to establish; a formal proof can be found in \cite[Theorem~3.1]{BeHoSh2012} (for continued fractions $(b_0,b_1,\dotsc)$; the proof for continued fractions $[b_0,b_1,\dotsc]$ is similar). Because of this correspondence between convergents and paths, we denote the sequence of convergents $v_0,v_1,\dotsc$ of a continued fraction by $\langle v_0,v_1,\dotsc\rangle $ (omitting $\infty$ as an initial vertex). 

This connection between continued fractions and paths in the Farey graph has been explored before, in, for example, \cite{BeHoSh2012,ShWa2016}. Using this perspective we can describe the function $\Phi$ in geometric terms, by looking at how it modifies the sequence of convergents. 

To see this, suppose first that $b_m=0$. Then $s_m(z)=-1/z$, so 
\[
S_m(\infty)=S_{m-1}\circ s_m(\infty)=S_{m-1}(0),
\]
and hence 
\[
v_m=S_m(\infty)=S_{m-1}(0)=S_{m-2}(\infty)=v_{m-2}.
\]
In this case the path of convergents $\langle v_0,v_1,\dotsc\rangle$ travels from $v_{m-2}$ to $v_{m-1}$ and then back to $v_{m-2}=v_m$, as shown on a schematic representation of part of $\mathscr{F}$ in the left-hand diagram of  Figure~\ref{fig2}. The effect of applying $\Phi$ is to remove the vertices $v_{m-1}$ and $v_m$ from the path of convergents. This will be proved formally in Lemma~\ref{lem5}.

\tikzset{every node/.style={circle,fill=black,draw=none,inner sep=2pt,label distance=-7pt}}
\tikzset{>=latex}

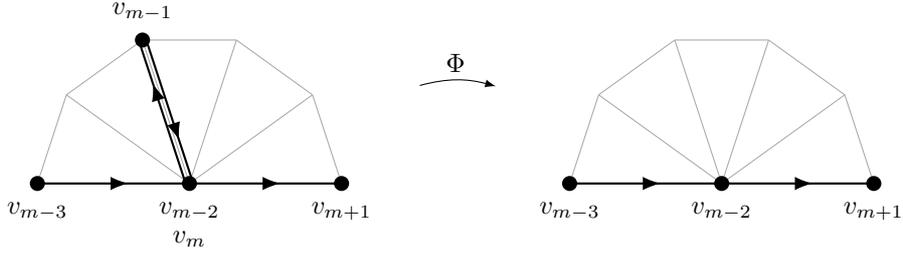
\begin{figure}[ht]
\centering
\begin{tikzpicture}[scale=1]
\pgfmathsetmacro{\ra}{2}
\pgfmathsetmacro{\tr}{5}
\pgfmathsetmacro{\eps}{0.05}

\begin{scope}[xshift=-3.5cm]
\foreach \n in {0,1,...,\tr} {
	\draw[grey] (0,0) -- (\n*180/\tr:\ra);
	\ifnum \n<\tr { 
		\draw[grey] (\n*180/\tr:\ra) --(\n*180/\tr+180/\tr:\ra);
	}
	\fi
}

\node[label={[]below:$v_{m-3}$}] at (-\ra,0){};
\node[label={[label={[yshift=7pt]below:$v_m$}]below:$v_{m-2}$}] at (0,0){};
\node[label={[]above:$v_{m-1}$}] at (3*180/\tr:\ra){};
\node[label={[]below:$v_{m+1}$}] at (\ra,0){};

\draw[thick,directed] (-\ra,0) -- (0,0);
\draw[thick,directedplus] ($(0,0)-(\eps,0)$) -- ($(3*180/\tr:\ra)-(\eps,0)$) ;
\draw[thick,directedplus]  ($(3*180/\tr:\ra)+(\eps,0)$)  -- ($(0,0)+(\eps,0)$) ;
\draw[thick,directed] (0,0) -- (\ra,0);
\end{scope}

\draw[->] (-0.47,1.3) to [out = 15, in = 165] (0.53,1.29);

\node[draw=none,fill=none] at (0,1.6){$\Phi$};

\begin{scope}[xshift=3.5cm]
\foreach \n in {0,1,...,\tr} {
	\draw[grey] (0,0) -- (\n*180/\tr:\ra);
	\ifnum \n<\tr { 
		\draw[grey] (\n*180/\tr:\ra) --(\n*180/\tr+180/\tr:\ra);
	}
	\fi
}

\node[label={[]below:$v_{m-3}$}] at (-\ra,0){};
\node[label={[]below:$v_{m-2}$}] at (0,0){};
\node[label={[]below:$v_{m+1}$}] at (\ra,0){};

\draw[thick,directed] (-\ra,0) -- (0,0);
\draw[thick,directed] (0,0) -- (\ra,0);
\end{scope}

\end{tikzpicture}
\caption{Geometric effect of the function $\Phi$ when $b_m=0$}
\label{fig2}
\end{figure}

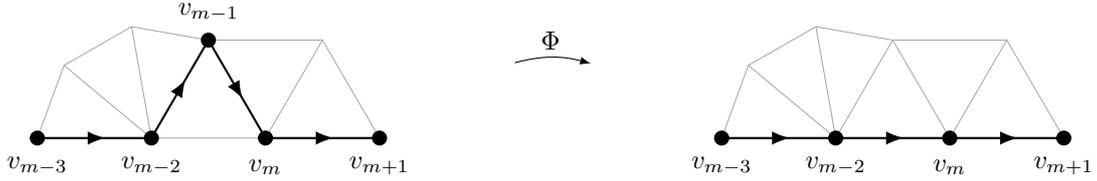
\begin{figure}[ht]
\centering
\begin{tikzpicture}[scale=1]
\pgfmathsetmacro{\ra}{1.5}
\pgfmathsetmacro{\eps}{0.05}
\pgfmathsetmacro{\h}{0.866*\ra}
\pgfmathsetmacro{\s}{0.684*\ra}

\begin{scope}[xshift=-4.5cm]

\draw[grey] (-1.5*\ra,-\h) -- (1.5*\ra,-\h) -- (\ra,0) -- (0.5*\ra,-\h) -- (0,0) -- (\ra,0);
\draw[grey] (-0.5*\ra,-\h) --++ (-220:\ra) --(-1.5*\ra,-\h);
\draw[grey] (-0.5*\ra,-\h) --++ (-260:\ra)--++(-150:\s);
\draw[grey] (-0.5*\ra,-\h) --++ (-300:\ra)--++(-190:\s);

\node[label={below:$v_{m-3}$}] at (-1.5*\ra,-\h){};
\node[label={below:$v_{m-2}$}] at (-0.5*\ra,-\h){};
\node[label={above:$v_{m-1}$}] at (0,0){};
\node[label={below:$\phantom{x}v_{m}\phantom{x}$}] at (0.5*\ra,-\h){};
\node[label={below:$v_{m+1}$}] at (1.5*\ra,-\h){};

\draw[thick,directed] (-1.5*\ra,-\h) -- (-0.5*\ra,-\h);
\draw[thick,directed] (-0.5*\ra,-\h) -- (0,0);
\draw[thick,directed]  (0,0) -- (0.5*\ra,-\h);
\draw[thick,directed]   (0.5*\ra,-\h)-- (1.5*\ra,-\h);
\end{scope}

\draw[->] (-0.47,-0.3) to [out = 15, in = 165] (0.53,-0.31);

\node[draw=none,fill=none] at (0,0){$\Phi$};

\begin{scope}[xshift=4.5cm]
\draw[grey] (-1.5*\ra,-\h) -- (1.5*\ra,-\h) -- (\ra,0) -- (0.5*\ra,-\h) -- (0,0) -- (\ra,0);
\draw[grey] (-0.5*\ra,-\h) --++ (-220:\ra) --(-1.5*\ra,-\h);
\draw[grey] (-0.5*\ra,-\h) --++ (-260:\ra)--++(-150:\s);
\draw[grey] (-0.5*\ra,-\h) --++ (-300:\ra)--++(-190:\s);

\node[label={below:$v_{m-3}$}] at (-1.5*\ra,-\h){};
\node[label={below:$v_{m-2}$}] at (-0.5*\ra,-\h){};
\node[label={below:$\phantom{x}v_{m}\phantom{x}$}] at (0.5*\ra,-\h){};
\node[label={below:$v_{m+1}$}] at (1.5*\ra,-\h){};

\draw[thick,directed] (-1.5*\ra,-\h) -- (-0.5*\ra,-\h);
\draw[thick,directed]  (-0.5*\ra,-\h) -- (0.5*\ra,-\h);
\draw[thick,directed]   (0.5*\ra,-\h)-- (1.5*\ra,-\h);
\end{scope}

\end{tikzpicture}
\caption{Geometric effect of the function $\Phi$ when $b_m=1$}
\label{fig3}
\end{figure}

Now suppose that $b_m=1$, in which case $s_m(1)=0$. Then $v_m=S_m(\infty)$, $v_{m-1}=S_{m-1}(\infty)=S_m(0)$ and 
\[
v_{m-2} = S_{m-2}(\infty)=S_{m-1}(0)=S_m(1).
\]
Since $1$, $0$ and $\infty$ are the vertices of a triangle in $\mathscr{F}$, so are $v_{m-2}$, $v_{m-1}$ and $v_m$. The effect of applying $\Phi$ is to skip the vertex $v_{m-1}$ and proceed directly from $v_{m-2}$ to $v_m$, as shown in Figure~\ref{fig3}. Again, this will be proved formally in Lemma~\ref{lem5} (and the case $b_m=-1$ is similar, but with a left rather than a right turn at $v_{m-1}$).

We also use the geometry of $\mathscr{F}$ to establish the following theorem about convergence of paths in $\mathscr{F}$. It is a generalisation of  \cite[Theorem~1.4]{ShWa2016}, which states (in a more general context) that if an infinite path in $\mathscr{F}$ does not return to any vertex infinitely often, then it converges. 

\begin{theorem}\label{thm3}
An infinite path in the Farey graph converges if and only if it does not return to any two distinct vertices infinitely often.
\end{theorem}

In the language of continued fractions, this theorem says that a continued fraction $[b_0,b_1,\dotsc]$ converges if and only if there are not two distinct extended rationals that each appear infinitely many times in the sequence of convergents of $[b_0,b_1,\dotsc]$.

\section{Paths in the Farey graph} 

We prove the three results stated in the introduction in reverse order, beginning with Theorem~\ref{thm3}. The next lemma will be used several times.

\begin{lemma}\label{lem1}
Let $u$ and $v$ be two adjacent vertices of $\mathscr{F}$, and let $\gamma$ be a finite path in $\mathscr{F}$ with initial and final vertices in different components of $\mathbb{R}_\infty\setminus \{u,v\}$. Then $\gamma$ passes through one or both of $u$ and $v$. 
\end{lemma}
\begin{proof}
Observe that the two vertices of an edge of $\gamma$ cannot lie in different components of $\mathbb{R}_\infty \setminus \{u,v\}$, for if they did then this edge would intersect the edge of $\mathscr{F}$ between $u$ and $v$. It follows, then, that $\gamma$ must pass through $u$ or $v$.
\end{proof}

It is helpful to highlight another elementary lemma.

\begin{lemma}\label{lem2}
Let $\alpha$ and $\beta$ be distinct elements of $\mathbb{R}_\infty$ that are not adjacent vertices of $\mathscr{F}$. Then the hyperbolic line between $\alpha$ and $\beta$ intersects some edge of $\mathscr{F}$.
\end{lemma}
\begin{proof}
Let $\ell$ be the hyperbolic line between $\alpha$ and $\beta$. If $\alpha$ and $\beta$ are irrational, then $\ell$ must intersect some edge of $\mathscr{F}$, for otherwise the vertices of $\mathscr{F}$ in one of the components of $\mathbb{R}_\infty \setminus \{\alpha,\beta\}$ are disconnected from the vertices in the other component. On the other hand, if one of the two vertices ($\alpha$, say) is rational, then after applying an element of $\Gamma$ we can assume that $\alpha=\infty$, in which case $\ell$ intersects the edge between $n$ and $n+1$, where $n$ is the integer part of $\beta$.
\end{proof}

We can now prove Theorem~\ref{thm3}. 

\begin{proof}[Proof of Theorem~\ref{thm3}]
We will prove the contrapositive statement of Theorem~\ref{thm3}, that an infinite path in the Farey graph diverges if and only if it returns to two distinct vertices infinitely often.

Let $\gamma=\langle v_0,v_1,\dotsc\rangle$ be a path in $\mathscr{F}$ that diverges. Then $\gamma$ must have two convergent subsequences with distinct limit points $\alpha$ and $\beta$ in $\mathbb{R}_\infty$. 

Suppose for the moment that $\alpha$ and $\beta$ are not adjacent vertices of $\mathscr{F}$. By Lemma~\ref{lem2}, there is an edge of $\mathscr{F}$ that intersects the hyperbolic line between $\alpha$ and $\beta$. Let $u$ and $v$ be the vertices of this edge. Then $\alpha$ and $\beta$ lie in different components of $\mathbb{R}_\infty\setminus \{u,v\}$. Since $\gamma$ approaches  each of $\alpha$ and $\beta$ infinitely often, we see from Lemma~\ref{lem1} that $\gamma$ passes through one of $u$ or $v$ infinitely many times.

Suppose now that $\alpha$ and $\beta$ are adjacent vertices of $\mathscr{F}$. After applying an element of the modular group $\Gamma$ we can assume that they are $0$ and $\infty$. Then $\gamma$ must pass in and out of one of the intervals $[-1,0]$ or $[0,1]$ infinitely often. Applying Lemma~\ref{lem1} once more, we see again that $\gamma$ passes through a vertex of $\mathscr{F}$ infinitely many times.

Thus, in both cases, and after applying another element of $\Gamma$, we can assume that $\gamma$ passes through the vertex $\infty$ infinitely often. However, $\gamma$ diverges, so it must enter some interval $[n,n+1]$ infinitely often, where $n$ is an integer. Since $n$ and $n+1$ are adjacent vertices of $\mathscr{F}$, we can apply Lemma~\ref{lem1} yet again, to see that $\gamma$ passes through one of $n$ or $n+1$ infinitely often. 

Therefore $\gamma$ returns to two distinct vertices of $\mathscr{F}$ infinitely often, as required. The converse implication is immediate.
\end{proof} 

\section{Convergent continued fractions}

Here we prove Theorem~\ref{thm2}. Although this theorem is a corollary of Theorem~\ref{thm1}, we prove it independently, and then later use it to prove the stronger theorem.

\begin{lemma}\label{lem3}
The continued fraction $[b_0,b_1,\dotsc]$, where $|b_n|\geq 2$, for $n=1,2,\dotsc$, converges to a value in $[b_0-1,b_0+1]$. Furthermore, it converges to $b_0-1$ if and only if  $b_1=b_2=\dotsb=2$, and it converges to $b_0+1$ if and only if $b_1=b_2=\dotsb=-2$.
\end{lemma}
\begin{proof}
We prove the lemma when $b_0=0$; the more general case follows immediately by applying a translation.

Let $t_n(z)=-1/(b_n+z)$, for $n=1,2,\dotsc$, and let $T_n=t_1\circ t_2\circ \dots \circ t_n$. (It is marginally more convenient to use these maps in place of the maps $s_n=b_n-1/z$ and $S_n=s_1\circ s_2\circ \dots \circ s_n$ defined in the introduction.) Then the $n$th convergent of $[0,b_1,b_2\dotsc]$ is $T_n(0)$, for $n=1,2,\dotsc$. Observe that $t_n$  maps the interval $[-1,1]$ inside itself and it preserves the order of points in that interval. It follows that $T_n$ also maps $[-1,1]$ inside itself and preserves order in that interval. Furthermore, $t_n(-1)=-1$ if and only if $b_n=2$, and $t_n(1)=1$ if and only if $b_n=-2$. Thus
\[
-1\leq T_1(-1)\leq T_2(-1) \leq \dotsb  \leq T_2(1)\leq T_1(1)\leq 1,
\]
where equality holds in all of the left set of inequalities if and only if all coefficients $b_n$ equal $2$, and equality holds in all of the right set of inequalities if and only if all coefficients $b_n$ equal $-2$. 

Now, $T_n(0)\in (T_n(-1),T_n(1))$, but $T_n(0)=T_{n+1}(\infty)$, so $T_{n}(0)\notin [T_{n+1}(-1),T_{n+1}(1)]$. Therefore either $T_n(-1)<T_n(0)<T_{n+1}(-1)$ or $T_{n+1}(1)<T_n(0)<T_n(1)$. From this we see that all the points $T_1(0),T_2(0), \dotsc$ are distinct. For any $\varepsilon>0$ there are only finitely many edges of $\mathscr{F}$ with vertices in $[-1,1]$ and with Euclidean diameter greater than $\varepsilon$. Since $T_n(0),T_n(1)\in [-1,1]$, and they are the vertices of an edge of $\mathscr{F}$ (because $0$ and $1$ are adjacent in $\mathscr{F}$ and $T_n\in \Gamma$), we see that $|T_n(0)-T_n(1)|\to 0$, and similarly $|T_n(0)-T_n(-1)|\to 0$. 

Reasoning in this way we deduce that the sequences $(T_n(-1))$ and $(T_n(1))$ converge to the same limit between $-1$ and $1$, which is the value of the continued fraction. Furthermore, this limit is $-1$ if and only if all coefficients $b_n$ equal $2$, and it is $1$ if and only if all coefficients $b_n$ equal $-2$. 
\end{proof}

We can now prove Theorem~\ref{thm2}, which says that if $|b_n|\geq 2$, for $n=1,2,\dotsc$, then the continued fraction $[b_0,b_1,\dotsc]$ converges to a rational if it has the same tail as $\pm[2,2,\dotsc]$, and otherwise it converges to an irrational.

\begin{proof}[Proof of Theorem~\ref{thm2}]
We use the notation $s_n(z)=b_n-1/z$, $S_n=s_1\circ s_2\circ \dotsb \circ s_n$ and $v_n=S_n(\infty)$, for $n=0,1,\dotsc$, which was presented in the introduction. The  convergents $\langle v_0,v_1,\dotsc\rangle$ of the continued fraction form a path in $\mathscr{F}$, which converges to a limit $\alpha$, a real number, by Lemma~\ref{lem3}.

Suppose that $\alpha$ is rational. We claim that there is a nonnegative integer $m$ for which $v_m\sim\alpha$ and $v_{m+1}\neq\alpha$. (Recall that $\sim$ denotes adjacency in $\mathscr{F}$.)  To prove the claim, first suppose that $v_k=\alpha$ for some nonnegative integer $k$. Then $v_{k+1}\sim\alpha$, and since $b_{k+2}\neq 0$ we see that $s_{k+2}(\infty)\neq s_{k+1}^{-1}(\infty)$, so
\[
v_{k+2} =S_{k+2}(\infty) = S_{k+1}\circ s_{k+2}(\infty) \neq S_{k+1}\circ s_{k+1}^{-1}(\infty) =S_k(\infty)=v_k = \alpha.
\]
Hence in this case we can choose $m=k+1$.

We can now suppose that no vertex $v_k$ is equal to $\alpha$. Choose vertices $u$ and $v$ of $\mathscr{F}$ each adjacent to $\alpha$ with $u<\alpha<v$ such that $v_0$ does not lie in $[u,v]$. By Lemma~\ref{lem1}, the path $\gamma$ must pass through one of $u$ or $v$, so there is a nonnegative integer $m$ for which $v_m\sim\alpha$ (and $v_{m+1}\neq \alpha$). This proves the claim.
 
Observe that $v_m=S_m(\infty)$, so $\infty=S_m^{-1}(v_m)$. Since $v_m\sim\alpha$, it follows that $\infty\sim S_m^{-1}(\alpha)$, so $S_m^{-1}(\alpha)$ is an integer. Clearly, $S_m^{-1}(\alpha)=[b_{m+1},b_{m+2},\dotsb]$. Next, $v_{m+1}=S_{m+1}(\infty)=S_m(b_{m+1})$, so $b_{m+1}=S_m^{-1}(v_{m+1})$, which is distinct from $S_m^{-1}(\alpha)$. We can now apply Lemma~\ref{lem3} to the continued fraction $[b_{m+1},b_{m+2},\dotsb]$ to see that either $S_m^{-1}(\alpha)=b_{m+1}-1$ and $b_{m+2}=b_{m+3}=\dotsb =2$, or else $S_m^{-1}(\alpha)=b_{m+1}+1$ and $b_{m+2}=b_{m+3}=\dotsb =-2$.

It remains to prove that if $[b_0,b_1,\dotsc]$ has the same tail as $[2,2,\dotsc]$ or $[-2,-2,\dotsc$], then $\alpha$ is rational. Suppose then  that $b_{m+2}=b_{m+3}=\dotsb =2$, for some nonnegative integer $m$. Then 
\[
\alpha = [b_0,b_1,\dotsc] = S_{m+1}([2,2,\dotsc])=S_{m+1}(1),
\]
so $\alpha$ is rational, and similarly we can see that $\alpha$ is  rational if $b_{m+2}=b_{m+3}=\dotsb =-2$.
\end{proof}

\section{Proof of the first part of Theorem~\ref{thm1}}

In this section we prove part~(i) of Theorem~\ref{thm1}. First we gather several more elementary results. In these results we use the usual notation $s_n(z)=b_n-1/z$ and $S_n=s_0\circ s_1\circ \dots \circ s_n$, for $n=0,1,\dotsc$, associated  to the continued fraction $[b_0,b_1,\dotsc]$.

\begin{lemma}\label{lem4}
Suppose that the finite continued fraction $[b_0,b_1,\dots,b_m]$, where $|b_i|\geq 2$ for $i=1,2,\dots,m$, has value $c/d$, a reduced rational. Then $m < d$.
\end{lemma}
\begin{proof}
From the equation $S_n=S_{n-1}\circ s_n$ a straightforward induction argument shows that $S_n(z) = (c_nz-c_{n-1})/(d_nz-d_{n-1})$ for coefficients $c_n$ and $d_n$ that satisfy 
\[
c_n=b_nc_{n-1}-c_{n-2},  \quad d_n=b_nd_{n-1}-d_{n-2} \quad\text{and}\quad c_{n-1}d_n-c_nd_{n-1}=1,
\]
for $n=1,2,\dots,m$, where $c_0=b_0$, $c_{-1}=1$, $d_0=1$ and $d_{-1}=0$. These are standard recurrence relations in continued fractions theory, perhaps most easily appreciated by representing the equation $S_n=S_{n-1}\circ s_n$ in matrix form as 
\[
\begin{pmatrix}c_n & -c_{n-1} \\ d_n & -d_{n-1} \end{pmatrix}=\begin{pmatrix}c_{n-1} & -c_{n-2} \\ d_{n-1} & -d_{n-2} \end{pmatrix}\begin{pmatrix}b_n & -1 \\ 1& 0 \end{pmatrix}.
\]
Using the recurrence relation for $d_n$ we see that
\[
|d_n| = |b_nd_{n-1}-d_{n-2}| \geq 2|d_{n-1}|-|d_{n-2}|,
\]
so $|d_n|-|d_{n-1}|\geq |d_{n-1}|-|d_{n-2}|$. Since $|d_0|-|d_{-1}|=1$, it follows that 
\[
|d_m| = (|d_m|-|d_{m-1}|)+(|d_{m-1}|-|d_{m-2}|)+\dots +(|d_0|-|d_{-1}|)>m.
\]
But $c_m/d_m=c/d$, so $|d_m|=d$, and the result follows. 
\end{proof}

The following lemma describes how the convergents of a continued fraction are modified under an application of the function $\Phi$.

\begin{lemma}\label{lem5}
Let $[b_0,b_1,\dotsc]$ and $[b_0',b_1',\dotsc]$ be continued fractions with convergents $\langle v_0,v_1,\dotsc \rangle$ and $\langle v_0',v_1',\dotsc\rangle$, respectively, and suppose that $[b_0',b_1',\dotsc]=\Phi([b_0,b_1,\dotsc])$. Let $m$ be the first positive integer position at which a $0$, $1$ or $-1$ appears in the sequence $b_1,b_2,\dotsc$ (assuming there is one). Then
\[
\langle v_0',v_1',\dotsc\rangle =
\begin{cases}
\langle v_0,v_1,\dots,v_{m-2},v_{m+1},v_{m+2},\dotsc\rangle, & \text{if $b_m=0$},\\
\langle v_0,v_1,\dots,v_{m-2},v_{m},v_{m+1},\dotsc\rangle, & \text{if $b_m=\pm1$}.
\end{cases}
\]
\end{lemma}

In the case $m=1$ this formula should be interpreted to say that $\langle v_0',v_1',\dotsc\rangle$ is equal to $\langle v_2,v_3,\dotsc\rangle$ when $b_m=0$ and $\langle v_0',v_1',\dotsc\rangle$ is equal to $\langle v_1,v_2,\dotsc\rangle$ when $b_m=\pm 1$.

\begin{proof}
We use the usual notation $s_n(z)=b_n-1/z$, $S_n=s_0\circ s_1\circ \dotsb \circ s_n$ and $v_n=S_n(\infty)$, for $n=0,1,\dotsc$. Also, we define $s_n'(z)=b_n'-1/z$, $S_n'=s_0'\circ s_1'\circ \dotsb \circ s_n'$ and $v_n'=S_n'(\infty)$, for $n=0,1,\dotsc$.

First consider the case when $b_m=0$.  Then $s'_n(z)=s_n(z)$, for $n=0,1,\dots,m-2$, so $v_n'=S_n'(\infty)=S_n(\infty)=v_n$. Now, $s_m(z)=-1/z$. Hence 
\[
s_{m-1}\circ s_m\circ s_{m+1}(z) = b_{m-1}+b_{m+1}-1/z=s_{m-1}'(z).
\]
 Also, $s_n'(z)=s_{n+2}(z)$, for $n\geq m$. Hence, for $n\geq m-1$, we have
\[
v_n' =S_n'(\infty) = S_{n+2}(\infty)=v_{n+2}.
\] 
Now suppose that $b_m=1$ (the case $b_m=-1$ is similar). Again, $s'_n(z)=s_n(z)$, for $n=0,1,\dots,m-2$, so  $v_n'=S_n'(\infty)=S_n(\infty)=v_n$. This time $s_m(z)=1-1/z$, and it can be checked that 
\[
s_{m-1}\circ s_m(\infty)=s_{m-1}'(\infty)\quad\text{and}\quad s_{m-1}\circ s_m\circ s_{m+1}(z) = s_{m-1}'\circ s_m'(z).
\]
Also, $s_n'(z)=s_{n+1}(z)$, for $n\geq m+1$.
 Hence, for $n\geq m-1$, we have
\[
v_n' =S_n'(\infty) = S_{n+1}(\infty)=v_{n+1}.\qedhere
\] 
\end{proof}

For the remainder of this section we will use the following notation associated to a continued fraction $[b_0,b_1,\dotsc]$. As usual, the sequence of convergents of $[b_0,b_1,\dotsc]$ is denoted by $\langle v_0,v_1,\dotsc\rangle$. We define $[b_0^{(n)},b_1^{(n)},\dotsc]=\Phi^n([b_0,b_1,\dotsc])$, and we let $\langle v_0^{(n)},v_1^{(n)},\dotsc\rangle$ be the associated sequence of convergents. Next we define $p^{(n)}$ to be the first positive integer position at which a 0, 1 or $-1$ appears in $[b_0^{(n)},b_1^{(n)},\dotsc]$ (the same definition appeared in the introduction), and we let $p=\liminf p^{(n)}$. By the definition of $\Phi$ we can see that, for $0\leq i<p-1$ (where possibly $p=\infty$), each of the sequences $b_i^{(0)},b_i^{(1)},\dotsc$ is eventually constant, with value $b_i^*$, say. We denote the sequence of convergents of $[b_0^*,b_1^*,\dots,b_{p-2}^*]$ by $\langle v_0^*,v_1^*,\dots,v_{p-2}^*\rangle$. Of course, if $p=\infty$, then $[b_0^*, b_1^*,\dotsc]$ is an infinite continued fraction with an infinite sequence of convergents $\langle v_0^*,v_1^*,\dotsc\rangle$. 

Consider now a particular convergent $v_k$ of $[b_0,b_1,\dotsc]$. We will define inductively a sequence $e^{(0)}(k),e^{(1)}(k),\dotsc$ of nonnegative integers chosen such that the vertex at position $e^{(n)}(k)$ of $\langle v_0^{(n)},v_1^{(n)},\dotsc\rangle$ is $v_k$. The resulting sequence may be finite or infinite.

First, let $e^{(0)}(k)=k$ (because $v_k$ is at position $k$ in $\langle v_0,v_1,\dotsc\rangle$). Now suppose that $e^{(0)}(k),e^{(1)}(k),\dots,e^{(n)}(k)$ have all been defined, for some nonnegative integer $n$. Let $m=p^{(n)}$, the position of the first coefficient (ignoring $b^{(n)}_0$) equal to $0$, $1$ or $-1$ in $[b_0^{(n)},b_1^{(n)},\dotsc]$. If $b_m^{(n)}=0$, and $e^{(n)}(k)$ equals $m-1$ or $m$, then the sequence terminates at $e^{(n)}(k)$. If $b_m^{(n)}=\pm 1$, and $e^{(n)}(k)=m-1$, then, again, the sequence terminates at $e^{(n)}(k)$. Otherwise, we define
\[
e^{(n+1)}(k) =
\begin{cases}
e^{(n)}(k), & \text{if $e^{(n)}(k)\leq m-2$,}\\
e^{(n)}(k)-2, & \text{if $e^{(n)}(k)> m$ and $b_m^{(n)}=0$,}\\
e^{(n)}(k)-1, & \text{if $e^{(n)}(k)> m-1$ and $b_m^{(n)}=\pm 1$.}\\
\end{cases}
\]
This definition has been chosen to ensure that the $e^{(n+1)}(k)$-th vertex of $\langle v_0^{(n+1)},v_1^{(n+1)},\dotsc\rangle$ is equal to the $e^{(n)}(k)$-th vertex of $\langle v_0^{(n)},v_1^{(n)},\dotsc\rangle$; that is, $v^{(n+1)}_{e^{(n+1)}(k)}=v^{(n)}_{e^{(n)}(k)}$. We can easily verify this by going through the cases, making use of Lemma~\ref{lem5}. 

The resulting sequence $e^{(0)}(k),e^{(1)}(k),\dotsc$ is (not necessarily strictly) decreasing, and if it is infinite, then it must therefore eventually be constant. It has the property that if $k<l$, then $e^{(n)}(k)<e^{(n)}(l)$, for all positive integers $n$ for which both expressions exist. This property is obvious intuitively, and it can be proved quickly from the definition.  We record one more property in a lemma. In this lemma we write $\mathbb{N}_0$ for the set $\{0,1,2,\dotsc\}$.

\begin{lemma}\label{lem7a}
Let $n,r\in\mathbb{N}_0$. Then there exists a unique integer $k\in\mathbb{N}_0$ with $e^{(n)}(k)=r$.
\end{lemma}

\begin{proof}
We prove by induction on $n$ that for any $r\in\mathbb{N}_0$ there exists $k\in\mathbb{N}_0$ with $e^{(n)}(k)=r$. For $n=0$, we can choose $k=r$, because $e^{(0)}(r)=r$. Suppose next that the induction statement is true for all nonnegative integers up to and including $n$. Let $m=p^{(n)}$. Given $r\in\mathbb{N}_0$, we define
\[
s =
\begin{cases}
r, & \text{if $r\leq m-2$,}\\
r+2, & \text{if $r\geq m-1$ and $b_m^{(n)}=0$,}\\
r+1, & \text{if $r\geq m-1$ and $b_m^{(n)}=\pm 1$.}
\end{cases}
\]  
By the inductive hypothesis, we can choose $k\in\mathbb{N}_0$ such that $e^{(n)}(k)=s$. Then one can check that $e^{(n+1)}(k)=r$ for each of the three cases used to define $s$. This completes the inductive proof. 

For uniqueness, we observe that if $k<l$ then $e^{(n)}(k)<e^{(n)}(l)$, so $e^{(n)}(k)$ and $e^{(n)}(l)$ cannot both equal $r$.
\end{proof}

All these properties of the sequence $e^{(0)}(k),e^{(1)}(k),\dotsc$ are intuitively clear, including the statement of the lemma, even if the details of formal proofs are somewhat awkward. Henceforth we often use these properties without comment.

\begin{lemma}\label{lem7}
Suppose that $p^{(n)}\rightarrow\infty$. Then for any extended rational $v$ there are only finitely many convergents $v_k$ equal to $v$. 
\end{lemma}
\begin{proof}
We write $v$ as a reduced rational $c/d$. Choose a positive integer $R$ such that if $n\geq R$, then $p^{(n)}\geq d+2$. It follows that 
\[
[b_0^{(n)},b_1^{(n)},\dots,b_{d}^{(n)}]=[b_0^*,b_1^*,\dots,b_{d}^*],
\]
for $n\geq R$. From Lemma~\ref{lem7a} there is a nonnegative integer $S$ such that $e^{(R)}(S)=d$. Then $e^{(R)}(S)\leq p^{(n)}-2$, so $e^{(n)}(S)\leq p^{(n)}-2$, for $n\geq R$, and it follows from the definition of the sequence $e^{(0)}(S),e^{(1)}(S),\dotsc$ that $e^{(n)}(S)=d$, for all $n\geq R$.

Let us assume that $k>S$ and $v_k=v$. We will establish a contradiction, thereby proving that only finitely many convergents are equal to $v$.

Suppose  for the moment that the sequence $e^{(0)}(k),e^{(1)}(k),\dotsc$ is infinite, so it is eventually constant, with value $m$, say. Choose any sufficiently large integer $n\geq R$ for which $e^{(n)}(k)=m$ and $p^{(n)}>m$. Since $k>S$, it follows that $e^{(n)}(k)>e^{(n)}(S)$, so $m>d$. Now, from the definition of $e^{(0)}(k),e^{(1)}(k),\dotsc$ we know that the vertex at position $m=e^{(n)}(k)$ of $\langle v_0^{(n)},v_1^{(n)},\dotsc \rangle$ is equal to $v_k$; that is, $v^{(n)}_m=v_k=v=c/d$. However, this contradicts  Lemma~\ref{lem4}, because $|b_i^{(n)}|\geq 2$, for $i=1,2,\dots,m$ (since $p^{(n)}>m$), and $m>d$.

Suppose instead that $e^{(0)}(k),e^{(1)}(k),\dotsc$ is of finite length -- the final term is $e^{(n)}(k)$, say. Since $k>S$, it follows that $e^{(n)}(k)> e^{(n)}(S)=d$.

Let $m=p^{(n)}$. Because the final term of the sequence $e^{(0)}(k),e^{(1)}(k),\dotsc$ is $e^{(n)}(k)$, we see from the definition of that sequence that either $b_m^{(n)}=0$ and $e^{(n)}(k)$ equals $m-1$ or $m$, or $b_m^{(n)}=\pm 1$ and $e^{(n)}(k)=m-1$. In all cases we have that $m>d$. 

Suppose  that $b_m^{(n)}=0$. If $e^{(n)}(k)=m-1$, then the finite continued fraction $[b_0^{(n)},b_1^{(n)},\dots,b_{m-1}^{(n)}]$ has value $v^{(n)}_{m-1}=v_k=v$, which contradicts Lemma~\ref{lem4}, because $|b_i^{(n)}|\geq 2$, for $i=1,2,\dots,m-1$ (since $p^{(n)}=m$), and $m>d$. And if $e^{(n)}(k)=m$, then the continued fraction $[b_0^{(n)},b_1^{(n)},\dots,b_{m-2}^{(n)}]$ has value $v^{(n)}_{m-2}=v^{(n)}_{m}=v_k=v$, which again contradicts Lemma~\ref{lem4}, because $|b_i^{(n)}|\geq 2$, for $i=1,2,\dots,m-2$ (since $p^{(n)}=m$), and $m>d$. We obtain a similar contradiction when $b_m^{(n)}=\pm 1$.

It follows that there are no integers $k>S$ with $v_k=v$. Hence there are only finitely many convergents equal to $v$.
\end{proof}

\begin{lemma}\label{lem6}
Suppose that $p^{(n)}\rightarrow \infty$. Then the sequence of convergents of the limit continued fraction $[b_0^*, b_1^*,\dotsc]$ is a subsequence of the sequence of convergents of $[b_0,b_1,\dots,b_n]$.
\end{lemma}
\begin{proof}
With the usual notation, the lemma says that $\langle v_0^*,v_1^*,\dotsc\rangle$ is a subsequence of $\langle v_0,v_1,\dotsc\rangle$. To see why this is so, choose any nonnegative integer $r$, and let $N$ be a positive integer for which $p^{(n)}\geq r+2$, for $n\geq N$. Then $v_r^{(n)}=v_r^*$, for $n\geq N$, because $r\leq p^{(n)}-2$. By Lemma~\ref{lem7a} there is a nonnegative integer $k_r$ with $e^{(N)}(k_r)=r$; and in fact $e^{(n)}(k_r)=r$, for $n\geq N$, because $r\leq p^{(n)}-2$. Hence $v_r^*=v_r^{(n)}=v_{k_r}$, by definition of $e^{(0)}(k_r),e^{(1)}(k_r),\dotsc$. 

Now, we know that if $k_r<k_s$, then $e^{(n)}(k_r)<e^{(n)}(k_s)$. Hence $k_0<k_1<\dotsb$. It follows that $\langle v_0^*,v_1^*,\dotsc\rangle$ is a subsequence of $\langle v_0,v_1,\dotsc\rangle$.
\end{proof}

We can now prove the first part of Theorem~\ref{thm1}.

\begin{proof}[Proof of Theorem~\ref{thm1}(i)]
Part (i) of Theorem~\ref{thm1} assumes that $p^{(n)}\to \infty$. In this case we see from Lemma~\ref{lem7} that the sequence of convergents $\langle v_0,v_1,\dotsc\rangle$ does not return to any vertex infinitely often. Hence, by Theorem~\ref{thm3}, the continued fraction converges to some value $\alpha$.

Now, by Lemma~\ref{lem6}, the sequence of convergents of $[b_0^*, b_1^*,\dotsc]$ is a subsequence of that of $[b_0,b_1,\dots]$. The former continued fraction certainly converges -- by Theorem~\ref{thm2} -- so it must converge to $\alpha$. Moreover, Theorem~\ref{thm2} tells us that $\alpha$ is rational if $[b_0^*, b_1^*,\dotsc]$ has the same tail as $[2,2,\dotsc]$ or $[-2,-2,\dotsc]$, and otherwise it is irrational, as required.
\end{proof}

\section{Proof of the second part of Theorem~\ref{thm1}}

The second part of Theorem~\ref{thm1} will be proved after two preliminary lemmas.

\begin{lemma}\label{lem8}
Consider a continued fraction $[b_0,b_1,\dotsc]$ with $|b_i|\geq 2$, for $i=1,2,\dots,n$. Let $\langle v_0,v_1,\dotsc\rangle$ be the corresponding sequence of convergents. If $1\leq k\leq n$, then 
\[
|v_{n}-v_{k-1}|\leq \frac{1}{|b_k|-1}.
\]
This inequality remains true if $1\leq k<n$ and $|b_i|\geq 2$, for $i=1,2,\dots,n-1$, but $b_n=0$.
\end{lemma}
\begin{proof}
By applying a translation we can assume that $b_0=0$ (so $v_0=0$). Now define $t_m(z)=-1/(b_m+z)$ and $T_m=t_1\circ t_2\circ \dotsb \circ t_m$, for $m=1,2,\dotsc$, as we did in proving Lemma~\ref{lem3}. Let $T_0$ denote the identity transformation. Recall that $(T_m(0))$ is the sequence of convergents of $[b_0,b_1,\dotsc]$ and $t_m([-1,1])\subset [-1,1]$, for $m=1,2,\dots,n$. Observe also that $|t_m'(z)|=1/|b_m+z|^2\leq 1$, for $z\in[-1,1]$ and $1\leq m\leq n$. Applying the chain rule, we see that
\[
|T_m'(z)| = |t_1'(z_1)||t_2'(z_2)|\dotsb |t_m'(z_m)|\leq 1,
\]
for $z\in[-1,1]$ and $1\leq m\leq n$, where $z_i=t_{i+1}\circ t_{i+2}\circ \dotsb \circ t_m(z)$ (and $z_m=z$).

Next, we have $v_k=T_k(0)$, for $k=1,2,\dots,n$. Observe that $t_{k+1}\circ t_{k+2}\circ \dotsb \circ t_n(0)\in [-1,1]$. Hence
\[
v_n \in T_k([-1,1]) = T_{k-1}([t_k(-1),t_k(1)]).
\]
Suppose that $b_k>0$ (the case $b_k<0$ is similar). Then $t_k(1)=-1/(b_k+1)<0$, so 
\[
v_n \in T_{k-1}([t_k(-1),0]).
\]
Therefore $v_n=T_{k-1}(u_n)$, for some point $u_n$ in $[t_k(-1),0]$. Now,
\[
|t_k(-1)-0| = \left| \frac{-1}{b_k-1}\right| = \frac{1}{|b_k|-1}.
\]
Hence
\[
|v_n-v_{k-1}|=|T_{k-1}(u_n)-T_{k-1}(0)|=|T_{k-1}'(c_k)||u_n-0|,
\]
for some real number $c_k$ between $u_n$ and $0$. The required inequality follows, since $|T_{k-1}'(c_k)|\leq 1$ and $|u_n-0|\leq 1/(|b_k|-1)$.

It remains to prove the final part of the lemma in which $b_n=0$. In this case, $v_n=v_{n-2}$, so the lemma continues to hold if $k\leq n-2$. And clearly it also holds if $k=n-1$.\end{proof}

\begin{lemma}\label{lem9}
Suppose that $p^{(n)}\not\rightarrow\infty$. Then only finitely many of the sequences  $e^{(0)}(k),e^{(1)}(k),\dotsc$, for $k\in\mathbb{N}_0$, are of infinite length.
\end{lemma}
\begin{proof}
Let $p=\liminf p^{(n)}$. Suppose that the sequence $e^{(0)}(k),e^{(1)}(k),\dotsc$ is of infinite length, for some positive integer $k$. Choose a positive integer $N$ such that, for $n\geq N$, all terms $e^{(n)}(k)$ of the sequence are equal to some value $e$.

Now suppose, in order to reach a contradiction, that $e>p$. Choose an integer $n\geq N$ for which $p^{(n)}=p$. Then $b^{(n)}_p$ is $0$, $1$ or $-1$, and we see from the definition of the sequence $e^{(0)}(k),e^{(1)}(k),\dotsc$  that $e^{(n+1)}(k)<e^{(n)}(k)$. This is the contradiction we need, since $e^{(n+1)}(k)=e^{(n)}(k)=e$.

It follows, then, that $e\leq p$. But for each nonnegative integer $e\leq p$ there is a \emph{unique} integer $k$ such that $e^{(n)}(k)\rightarrow e$ as $n\rightarrow\infty$. Hence there are only finitely many infinite sequences $e^{(0)}(k),e^{(1)}(k),\dotsc$, for $k\in\mathbb{N}_0$.
\end{proof}

We now prove the second part of Theorem~\ref{thm1}.

\begin{proof}[Proof of Theorem~\ref{thm1}(ii)]
Part (ii) of Theorem~\ref{thm1} assumes that $p^{(n)}\not\rightarrow \infty$. Let $p=\liminf p^{(n)}$. 

Suppose that $q^{(n)}\not\rightarrow\infty$. Since $p=\liminf p^{(n)}$, we can find a positive integer $N$ such that 
\[
[b_0^{(n)},b_1^{(n)},\dots,b_{p-2}^{(n)}]=[b_0^*,b_1^*,\dots,b_{p-2}^*],
\]
for $n\geq N$. Observe that the sequence $b_{p-1}^{(1)},b_{p-1}^{(2)},\dotsc$ does not stabilise on a fixed value. And recall that $q^{(n)}=|b_{p-1}^{(n)}|$, by definition.

Now let $m_1,m_2,\dotsc$ be the complete list of positive integers $n$ greater than $N$, written in increasing order, for which $p^{(n)}=p$. That is, $m_1,m_2,\dotsc$ are the positive integers $n>N$, in order, for which $b^{(n)}_p$ is $0$, $1$ or $-1$.

The sequence $b_{p-1}^{(n)}$, for $n=N+1,N+2,\dotsc$, can only change value when $n=m_1,m_2,\dotsc$. Since $q^{(n)}\not\rightarrow\infty$, we can find a subsequence $n_1,n_2,\dotsc$ of $m_1,m_2,\dotsc$ for which every term $b_{p-1}^{(n_i)}$ is equal to some fixed integer $b$, where $|b|\geq 2$. And by restricting to a further subsequence, we can assume that all the terms $b_{p}^{(n_i)}$ are equal to precisely one of $0$, $1$ or $-1$.

Observe that $v_{p-1}^{(n_i)}=u$ and $v_{p}^{(n_i)}=v$, for $i=1,2,\dotsc$, where $u$ and $v$ are two fixed adjacent vertices of $\mathscr{F}$. Let $r_i$ and $s_i$ be nonnegative integers for which $e^{(n_i)}(r_i)=p-1$ and $e^{(n_i)}(s_i)=p$, in which case $v_{r_i}=u$ and $v_{s_i}=v$. Now, the sequence $e^{(0)}(r_i),e^{(1)}(r_i),\dotsc$ has length exactly $n_i+1$; the final term is $e^{(n_i)}(r_i)=p-1$ because $b_{p}^{(n_i)}$ is $0$, $1$ or $-1$. It follows that all the integers $r_i$ are distinct from one another, so there are infinitely many of them. With similar reasoning we can see that the collection of integers $s_i$ is infinite too.

We deduce that the sequence $(v_n)$ is equal to $u$ for infinitely many indices $n$, and likewise it is equal to $v$ for infinitely many indices $n$. It follows that $(v_n)$ diverges, so $[b_1,b_2,\dotsc]$ diverges.


Suppose now that $q^{(n)}\rightarrow\infty$. As before, we choose a positive integer $N$ for which $p^{(n)}\geq p$, for $n\geq N$, so 
\[
[b_0^{(n)},b_1^{(n)},\dots,b_{p-2}^{(n)}]=[b_0^*,b_1^*,\dots,b_{p-2}^*].
\]
Next, by Lemma~\ref{lem9}, we can choose a positive integer $M$ such that, for $k\geq M$, the sequence $e^{(0)}(k),e^{(1)}(k),\dotsc$ is of finite length. Let us also choose $M>2N+(p-1)$. 

Let $k\geq M$. By definition of the sequence $e^{(0)}(k),e^{(1)}(k),\dotsc$, we can see that 
\[
k-e^{(n)}(k) =  (k-e^{(1)}(k))+(e^{(1)}(k)-e^{(2)}(k))+\dots +(e^{(n-1)}(k)-e^{(n)}(k))\leq 2n.
\]
Suppose, in order to reach a contradiction, that $e^{(n)}(k)<p-1$, for some positive integer $n$. Then $k-(p-1)<2n$, and since $M>2N+(p-1)$ and $k\geq M$, we have that $n>N$. However, if $e^{(n)}(k)\leq p-2\leq p^{(n)}-2$, then $e^{(n+1)}(k)=e^{(n)}(k)$, so the sequence $e^{(0)}(k),e^{(1)}(k),\dotsc$ is infinite. This is the required contradiction. Therefore $e^{(n)}(k)\geq p-1$ for $k\geq M$ and $n\in\mathbb{N}_0$.

Next, for each integer $k>N$, define $r_k$ to be the positive integer such that $e^{(r_k)}(k)$ is the last term of the sequence $e^{(0)}(k),e^{(1)}(k),\dotsc$. Let $n\in\mathbb{N}$ and $m=p^{(n)}$. Then $r_k=n$ if and only if either $b_m=0$ and $e^{(n)}(k)$ equals $m-1$ or $m$, or $b_m=\pm 1$ and $e^{(n)}(k)=m-1$. Hence there are at most two values of $k$ for which $r_k=n$. Consequently, we deduce that $r_k\rightarrow\infty$ as $k\rightarrow\infty$.

The vertex at position $e^{(r_k)}(k)$ of $\langle v^{(r_k)}_0,v^{(r_k)}_1,\dotsc\rangle$ is equal to $v_k$, so we can apply Lemma~\ref{lem8} to the continued fraction $[b^{(r_k)}_0,b^{(r_k)}_1,\dotsc]$ to see that
\[
|v_k-v^{(r_k)}_{p-2}| \leq \frac{1}{|b_{p-1}^{(r_k)}|-1}=\frac{1}{q^{(r_k)}-1}.
\]
Observe that $v^{(r_k)}_{p-2}=v_{p-2}^*$, for sufficiently large values of $k$. Since $r_k\to\infty$ as $k\to\infty$, and hence $q^{(r_k)}\to\infty$ as $k\to\infty$, we see that $v_k\to v_{p-2}^*$ as $k\to\infty$. Therefore $[b_0,b_1,\dotsc]$ converges to the extended rational $v_{p-2}^*$.
\end{proof}



\end{document}